\documentclass [preprint,11pt,a4paper]{elsarticle}%[11pt,a4paper]{amsart}

\newcommand{\pn}{\{ P_n(x) \}_{n=0}^\infty }
\newcommand{\qn}{\{ Q_n(x) \}_{n=0}^\infty }
\newcommand{\rn}{\{ R_n(x) \}_{n=0}^\infty }
\newcommand{\tn}{\{ T_n(x) \}_{n=0}^\infty }
\newcommand{\sn}{\{ S_n^{(\lambda)}(x) \}_{n=0}^\infty }
\newcommand{\hl}{^{(\ld)}}
\newcommand{\bl}{_{\lambda}}
\newcommand{\si}{ \sigma}
\newcommand{\ld}{ \lambda}
\newcommand{\la}{ \langle}
\newcommand{\ra}{ \rangle}
\newcommand{\uu}{\{u_0,u_1\}}
\newcommand{\dproof}{\noindent {\em Proof. \quad}}{}

\usepackage{amsmath,amsthm,amssymb,amsfonts}

\theoremstyle{plain}
\newtheorem{theorem}{Theorem}[section]

\theoremstyle{definition}
\newtheorem{definition}{Definition}[section]
\newtheorem{example}{Example}[section]

\theoremstyle{remark}
\newtheorem{remark}{Remark}[section]

\swapnumbers

\journal{Indagationes Mathematicae}

\begin{document}

\begin{frontmatter}

\title{On an extension of generalized coherent pairs of orthogonal polynomials: the classical case}

\author[inst1]{J. H. Lee}

\affiliation[inst1]{organization={ARIST},%Department and Organization
            addressline={85 Wolpyeongbuk-ro}, 
            city={Seo-gu},
            postcode={35213}, 
            state={Daejeon},
            country={Korea}}

\author[inst2]{S. J. An}

\affiliation[inst2]{organization={Dept. of Intelligent System Engineering},%Department and Organization
            addressline={Cheju Halla University, 38 Halladaehak-ro}, 
            city={Jeju-si},
            postcode={63092}, 
            state={Jeju Special Self-Governing Province},
            country={Korea}}
            
\author[inst3]{H. Y. Lee}

\affiliation[inst3]{organization={Dept. of Mathematics},%Department and Organization
            addressline={University of Utah Asia Campus, 119-3 Songdo Moonwha-ro}, 
            city={Yeonsu-Gu},
            postcode={21985}, 
            state={Incheon},
            country={Korea}}

%\subjclass{Primary 33C45; 42C05}
%\keywords{Orthogonal Polynomials, Coherent Pairs, Sobolev Inner Product}

\vskip 0.5cm

\begin{abstract}
{ Given two quasi-definite moment functionals, the corresponding orthogonal polynomial systems satisfy an algebraic differential relation(called an extended coherent pair).
We study generalizing extended coherent pairs that unify extended coherent pairs and extended symmetric coherent pairs and find the related coefficients.
When one of the moment functionals is (strongly) classical, we find another orthogonal polynomial system to find three-term recurrence coefficients.
Moreover, we determine the companion moment functional as a rational modification of the classical one.}

\end{abstract}

\begin{keyword}
%% keywords here, in the form: keyword \sep keyword
Orthogonal Polynomials \sep Coherent Pairs \sep Sobolev Inner Product
%keyword one \sep keyword two
%% PACS codes here, in the form: \PACS code \sep code
\PACS Primary 33C45 \sep 42C05
%% MSC codes here, in the form: \MSC code \sep code
%% or \MSC[2008] code \sep code (2000 is the default)
%\MSC 0000 \sep 1111
\end{keyword}

%\maketitle
\end{frontmatter}

\section{Introduction}
\setcounter{equation}{0} \setcounter{theorem}{0}

Concerning polynomials orthogonal with respect to a
Sobolev inner product
\begin{equation}\label{eq11}
\phi\bl(f,g):=\la u_0, fg\ra+\ld\la u_1,f'g'\ra~(\ld\neq 0),
\end{equation}
where $u_0$ and $u_1$ are quasi-definite moment functionals, 
Iserles et al.\cite{Is} introduced the concept of coherency and symmetric coherency for the pair of moment functionals $\uu$.
$\uu$ is a coherent pair(resp. symmetrically coherent pair) if there are nonzero constants $\si_n$(resp. $\tau_n$) such that
$$\aligned R_n(x)&=\frac{1}{n+1}P'_{n+1}(x)-\frac{\sigma_n}{n}P'_n(x),~n\geq 1\\
\text{(resp., }
R_n(x)&=\frac{1}{n+1}P'_{n+1}(x)-\frac{\tau_{n-1}}{n-1}P'_{n-1}(x),~n\geq2), \endaligned$$ 
where $\pn$ and $\rn$ are the monic orthogonal
polynomials with respect to $u_0$ and $u_1$, respectively,
and in \cite{Me1}, all coherent pairs are classified.

After the work by Iserles et al.\cite{Is}, the coherency was
studied in various ways. One is concerned with generalizing the coherency. In \cite{KLM}, generalized coherency, which unifies both coherency and symmetric coherency, was studied.
Generalized coherency is given by the relation,
\begin{equation}\label{eq12}
R_n(x)=Q_n(x)-\sigma_nQ_{n-1}(x)-\tau_{n-1}Q_{n-2}(x),~~n\geq 0,
\end{equation}
where $Q_n(x)=\frac{1}{n+1}P'_{n+1}(x),~n\geq 0$, $Q_{-1}(x)=Q_{-2}(x)=0$, and $\sigma_0=\tau_{-1}=\tau_0=0$.

In this case, the monic orthogonal polynomials $\sn$ with respect to
$\phi\bl(\cdot,\cdot)$ and $\pn$ satisfy the relation
\begin{multline}\label{eq13}
P_{n+1}(x)-{\sigma}_nP_n(x)-{\tau}_{n-1}P_{n-1}(x)\\
=S_{n+1}^{(\ld)}(x)-\mu_nS_n^{(\ld)}(x)-\theta_{n-1}S_{n-1}^{(\ld)}(x),
~n\geq 0,
\end{multline}
where 
$\tau_{-1}= \theta_{-1}= P_{-1}(x)=S_{-1}\hl(x)=0$,
and $\mu_n$ and $\theta_n$ are evaluated.

Another approach for generalizing coherency has occurred. 
Extended coherency was introduced and studied  in \cite{KKMY} 
when $\pn$ and $\sn$ have the relation
\begin{equation}\label{eq14}
P_{n+1}(x)-\sigma_nP_n(x) =S_{n+1}^{(\ld)}(x)-\mu_nS_n^{(\ld)}(x), ~n\geq 0.
\end{equation}
Note that this relation is the 2-term case of (\ref{eq13}), and $\pn$ and $\rn$
satisfy the relation
$$ R_n(x)-d_{n}R_{n-1}(x) =\frac{1}{n+1}P'_{n+1}(x) -\frac{\sigma_n}{n}P'_n(x), ~n\geq 1,$$
for some $d_n$. We call this extended coherency.

In \cite{DM2}, extended symmetric coherency was studied.
When $\pn$ and $\sn$ have the relation
\begin{equation}\label{eq15}
P_{n+1}(x)-\tau_{n-1}P_{n-1}(x)
=S_{n+1}^{(\ld)}(x)-\theta_{n-1}S_{n-1}^{(\ld)}(x), ~n\geq 0,
\end{equation}
note that this relation is the symmetric case of (\ref{eq13}),
and then $\pn$ and $\rn$ satisfy the relation
$$ R_n(x)-e_{n-1}R_{n-2}(x) =\frac{1}{n+1}P'_{n+1}(x)-\frac{\tau_{n-1}}{n-1}P'_{n-1}(x),~n\geq 2,$$
for some $e_n$. We call this extended symmetric coherency.

\vskip 0.3cm
The concept of coherency is based on Sobolev orthogonality  and is applicable for Fourier-Sobolev expansions. This brought a lot of attention in the area of coherency and the research of extending coherency has been studied in many different directions. 

In \cite{BFMR,MS}, the concept of coherent pairs of measures was extended from the real line to Jordan arcs and curves,
focusing on the coherent pairs supported on the unit circle.

In \cite{JP, Pe}, $(M,N)-$coherency was introduced as extending terms of polynomials. 
$(M,N)-$coherency has the relation
$$ \sum_{i=0}^N r_{i,n} R_{n-i}(x) = \sum_{i=0}^M p_{i,n} P'_{n+1-i}, ~n\geq 0,$$
with the conditions $r_{N,n}\neq 0$ if $n \geq N$ and $p_{M,n}\neq 0$ if $n \geq M$, and the conventions
$r_{i,n} = 0$ if $n < i \leq N$ and $p_{i,n} = 0$ if $n < i \leq M$. 

In \cite{CM, JMPP, MP}, coherency was extended to the relation of higher-order differentials of the corresponding polynomials. Combining with $(M,N)-$coherency, $(M,N)-$coherency of order $(m,k)$ has the relation
$$ \sum_{i=0}^N r_{i,n} R^{(k)}_{n-i}(x) = \sum_{i=0}^M p_{i,n} P^{(m)}_{n+1-i}, ~n\geq 0,$$
where $r_{i,n}$ and $p_{j,n}$ are complex numbers with $r_{N,n}\neq 0$ if $n \geq N$ and $p_{M,n}\neq 0$ if $n \geq M$, and $r_{i,n} =p_{i,n}=0$ if $i>n$.

In \cite{AMPR}, linearly related orthogonal polynomials and their functionals are studied, and which can be applied to the linear relation of two polynomials.
In \cite{DM1}, companion of linear functional in Sobolev inner product is studied.

In \cite{MX}, Marcell\'an and Xu provided a survey of Sobolev orthogonal polynomials, including coherent pairs and generalizing coherent pairs.

\vskip 0.3cm

In this work, we will generalize extended coherency that unifies
extended coherency and extended symmetric coherency, that is, in the
case that $\sn$ and $\pn$ have the relation (\ref{eq13}), we have
the result that unifies the result of extended coherency and that of extended symmetric coherency.

In section 3, we introduce generalized extended coherency
and find the related coefficients.
In section 4, when one of the moment functionals is (strongly) classical, we
find the three-term recurrence relation coefficients of the monic
orthogonal polynomial system relative to the other moment functional. We find the companion moment functional if one is given.

\vskip .5cm

\section{Preliminaries}
\setcounter{equation}{0} \setcounter{theorem}{0}

Let $\mathbb{P}$ be the linear space of all polynomials in one
variable with complex coefficients. We denote the degree of a
polynomial $P(x)$ by $\deg(P)$ with the convention that
$\deg(0)=-1$. A polynomial system(PS) is a sequence of polynomials
$\pn$ with $\deg(P_n)=n,~n\geq0$. For convenience, let $P_n(x)=0$ if $n<0$, and we denote $P_n$ instead of $P_n(x)$ for simplicity. 

A linear functional $u$ on $\mathbb{P}$ is called a moment
functional, and we denote its action on a polynomial $\phi(x)$ by
$\la u,\phi\ra$. We say that a moment functional $u$ is
quasi-definite(positive-definite, respectively) if its moments
$a_n:=\la u,x^n\ra,~n\geq0$, satisfy the Hamburger condition
$$ \Delta_n(u):=\det[a_{i+j}]_{i,j=0}^{n}\neq 0, \quad (\Delta_n(u)>0,~\text{respectively}),\quad n\geq 0.$$

\begin{definition}\label{def2.1} 
A PS $\pn$ is said to be an orthogonal
polynomial system(OPS) if there is a linear functional $u$ on
$\mathbb{P}$ such that
$$\langle u,P_m P_n \rangle = p_n\delta_{mn}, ~ m,n\geq0,$$
where $p_n$ are non-zero constants. 

In this case, we call $\pn$ an
OPS relative to $u$, and $u$ is said to be an orthogonalizing
moment functional of $\pn$. 

A linear functional $u$ is
quasi-definite if and only if there is an OPS $\pn$ relative to
$u$ (see \cite{Ch2}). 

Moreover, in this case, each $P_n(x)$ is
uniquely determined up to a non-zero constant factor.

A polynomial $P_n(x)$ is monic if the leading coefficient of $P_n(x)$ is equal to $1$. 
We call $\pn$ an MOPS(monic OPS) if all polynomials in $\pn$ are monic.

\end{definition}

It is well-known(\cite{Ch2}) that such sequences of monic
polynomials satisfy three-term recurrence relations
\begin{equation}
P_{n+1}(x)=(x-b_n)P_n(x)-c_nP_{n-1}(x),~n\geq 0,
\end{equation}
where $P_{-1}(x)=0$, $P_0(x)=1$, $c_n\neq 0$ for $n\geq 0$, and
$c_0=\la u, 1\ra$.

For a moment functional $u$, a polynomial $\phi(x)$, and a
constant $c$, we define moment functionals $u',~\phi u~$, and
$(x-c)^{-1}u$ by
\begin{align*}
\langle u',p \rangle&=-\langle u,p' \rangle; \\
\langle \phi u,p \rangle&=\langle u,\phi p\rangle;\\
\langle (x-c)^{-1}u,p \rangle&=\langle
u,\frac{p(x)-p(c)}{x-c}\rangle,\quad p\in \mathbb{P}.
\end{align*}

\begin{definition}\label{def2.2}(\cite{Ma}) 
A quasi-definite moment functional $u$ is said to be semiclassical if $u$ satisfies
\begin{equation}\label{eq21}
(\varphi u)'=\psi u,
\end{equation}
for some polynomials $\varphi(x)$ and $\psi(x)$ with
$(\varphi,\psi)\neq (0,0)$. We then have $\deg(\varphi)\geq 0$ and
$\deg(\psi)\geq 1$. The corresponding OPS is called a
semiclassical OPS.
\end{definition}

For a semiclassical moment functional $u$,
$$s := \min\max (\deg(\varphi)-2, \deg(\psi)-1)$$
the class number of $u$, where the minimum is taken over all pairs
$(\varphi,\psi) \ne (0,0)$ of polynomials satisfying (\ref{eq21}).
In particular, a semiclassical moment functional of class $0$ is
called a classical moment functional.

It is well-known that there are essentially four distinct classical OPS's, up to a linear change of variable(\cite{Bo,KL}), 
for each case, we denote the corresponding orthogonalizing moment functional by $u$ :
\begin{itemize}
\item[(i)] Hermite polynomials $\{H_n(x)\}_{n=0}^\infty$: $\varphi(x)= 1$, $\psi(x)=-2x$, $u=u_H$;
\item[(ii)] Laguerre polynomials $\{L_n^{(\alpha)}(x)\}_{n=0}^\infty$: $\varphi(x)=x$, $\psi(x)=-x+\alpha+1~ (\alpha\notin \{-1,-2,\cdots\}$),$u=u_L^{(\alpha)}$;
\item[(iii)] Bessel polynomials $\{B_n^{(\alpha)}(x)\}_{n=0}^\infty$: $\varphi(x)=x^2$,
$\psi(x)=(\alpha+2)x+2$ ($\alpha\notin\{-2,-3,$ $\cdots\}$), $u=u_B^{(\alpha)}$;
\item[(iv)] Jacobi polynomials $\{P_n^{(\alpha,\beta)}(x)\}_{n=0}^\infty$: $\varphi(x)=1-x^2$,
$\psi(x)=\beta-\alpha-(\alpha+\beta+2)x$ ($\alpha,\beta,\alpha+\beta+1\notin\{-1,-2,\cdots\}$), $u=u_J^{(\alpha,\beta)}$.
\end{itemize}

We say a quasi-definite moment functional $u$ with MOPS $\pn$ to be strongly classical(see \cite{KLM}) if
there is another MOPS $\tn$ relative to $w$ such that $ P_n(x)= \frac{1}{n+1}T'_{n+1}(x),~n\geq 0$.
Then $u$ and $w$ must be classical moment functionals of the same type satisfying
$$ (\varphi u)'=\psi u,~~ (\varphi w)'=(\psi-\varphi')w,\text{ and } \varphi w=u.$$
Classical moment functionals $u_J^{(\alpha,\beta)}$ $(\alpha, \beta, \alpha+\beta \neq 0,-1,-2,\cdots)$,
$u_B^{(\alpha)}$ $(\alpha\neq0,-1,$ $ -2, \cdots)$, $u_L^{(\alpha)}~(\alpha\neq 0,-1,-2,\cdots)$, and $u_H$ are strongly classical.

\begin{definition}(\cite{KLM}) \label{def2.3} 
Let $u_0$ and $u_1$ be quasi-definite
moment functionals with corresponding MOPS $\pn$ and $\rn$,
respectively. Let $Q_n(x)=\frac{1}{n+1}P'_{n+1}(x),$ $n\geq 0.$ 
\begin{itemize}
\item[(i)] $\uu$ is a coherent pair if there exist
complex numbers $\sigma_n$, $n\geq1$, such that
$$ R_n(x)=Q_n(x) -\sigma_{n-1} Q_{n-1}(x),~~n\geq 0, $$
where $\sigma_{-1}=0$.

\item[(ii)] Assume $u_0$ and $u_1$ are symmetric moment
functionals. $\uu$ is a symmetrically coherent pair if there exist
complex numbers $\tau_n$, $n\geq1$, such that
$$ R_n(x)=Q_n(x)-\tau_{n-2}Q_{n-2}(x),~~n\geq 0, $$
where $\tau_{-1}=\tau_{-2}=0$.

\item[(iii)] $\uu$ is a generalized coherent pair if there
exist complex numbers $\sigma_n$ and
$\tau_n$, $n\geq1$, such that
$$ R_n(x)=Q_{n}(x)-\sigma_{n-1} Q_{n-1}(x)-\tau_{n-2}Q_{n-2}(x),~~n\geq 0, $$
where $\sigma_{-1}=\tau_{-1}=\tau_{-2}=0$.
\end{itemize}
\end{definition}

\begin{definition}(\cite{DM2, KKMY}) \label{def2.4} 
Let $u_0$ and $u_1$ be quasi-definite
moment functionals with corresponding MOPS $\pn$ and $\rn$,
respectively. Let $Q_n(x)=\frac{1}{n+1}P'_{n+1}(x),$ $n\geq 0.$
\begin{itemize}
\item[(i)] $\uu$ is an extended (2-term) coherent pair if there exist
complex numbers $d_n$ and $\sigma_n$, $n\geq0$, such that
$$ R_n(x)-d_{n-1}R_{n-1}(x) =Q_n(x) -\sigma_{n-1} Q_{n-1}(x),~~n\geq 0, $$
where $d_{-1}=\sigma_{-1}=0$.

\item[(ii)] Assume $u_0$ and $u_1$ are symmetric moment
functionals. $\uu$ is an extended symmetrically coherent pair if there exist complex numbers $e_n$ and $\tau_n$, $n\geq0$, such that
$$ R_n(x)-e_{n-2}R_{n-2}(x) =Q_n(x)-\tau_{n-2}Q_{n-2}(x),~~n\geq 0, $$
where $e_{-1}=e_{-2}=\tau_{-1}=\tau_{-2}=0$.
\end{itemize}
\end{definition}

\vskip .5cm

\section{Extended Generalized Coherent Pairs}
\setcounter{equation}{0} \setcounter{theorem}{0}

Let $u_0$ and $u_1$ be two quasi-definite moment functionals with
corresponding MOPS's $\pn$ and $\rn$, respectively. Let $\sn$ be
the MOPS relative to the Sobolev inner product (\ref{eq11}). We define
$$ p_n:=\langle u_0,P_{n}^{2} \rangle,~
r_n:=\langle u_1,R_{n}^{2} \rangle,~
s_n:=\phi_\ld(S\hl_n,S\hl_n),~n\geq 0.$$

When $\uu$ is a generalized coherent pair, $\pn$ and $\rn$ satisfy
the relation in Definition \ref{def2.3} (iii).

In this case, $\sn$ and $\pn$
satisfy the relation (see \cite{KLM}), for some $\tilde\sigma_n, \tilde\tau_n,\mu_n,$ and $\theta_n,$ $n\geq 0$
\begin{multline}\label{eq31}
P_{n+1}(x)-\tilde\sigma_nP_n(x)-\tilde\tau_{n-1}P_{n-1}(x) \\
=S_{n+1}^{(\ld)}(x)-\mu_nS_n^{(\ld)}(x)-\theta_{n-1}S_{n-1}^{(\ld)}(x),
~n\geq 0,
\end{multline}
where $\theta_{-1}=\tilde\tau_{-1}=0$ and $\mu_n$ and $\theta_n$ satisfy the relations
\begin{equation} \label{eq32}
\left\{ \begin{aligned}
\mu_n s_n&=\tilde\sigma_n p_n + \tilde\tau_{n-1}(\mu_{n-1}-\tilde\sigma_{n-1})
p_{n-1}, ~n\geq 1, \\
\theta_ns_n &=\tilde\tau_n p_n, ~n\geq 1.
\end{aligned}\right.
\end{equation}
Now we consider the inverse case. We assume that $\sn$
and $\pn$ satisfy the relation (\ref{eq31}).
\begin{theorem}\label{thm31}
Assume $\pn$ and $\sn$ satisfy the relation (\ref{eq31}), then $\pn$ and $\rn$ satisfy the relation \begin{multline}\label{eq33}
(n+1)R_n(x)-\tilde d_{n}R_{n-1}(x)-\tilde e_{n-1}R_{n-2}(x) \\
=P'_{n+1}(x)-\tilde\sigma_n P'_n(x)-\tilde\tau_{n-1}P'_{n-1}(x),~n\geq 0,
\end{multline}
for some $\tilde d_n$ and $\tilde e_n$, $n\geq 1$ with $\tilde d_0=\tilde e_0=\tilde e_{-1}=0$ and 
\begin{align}
&\begin{aligned} 
\tilde d_{n}=&\frac{\mu_n s_n-\tilde\sigma_np_n
+(\tilde\sigma_{n-1} -\mu_{n-1})\tilde\tau_{n-1}p_{n-1}}{n\ld r_{n-1}} \\
&+\frac{r_{n-2}\tilde e_{n-1}\{\tilde d_{n-1}-(n-1)\mu_{n-1}\}}{nr_{n-1}}
,~n\geq 2, \\
\tilde d_1= &\frac{\mu_1 s_1-\tilde\sigma_1 p_1 +(\tilde\sigma_{0} -\mu_{0})\tilde\tau_{0}p_{0}}{\ld r_{0}}, \label{eq34} \end{aligned} \\
& \tilde e_{n-1}=\frac{\theta_{n-1}s_{n-1}-\tilde\tau_{n-1}p_{n-1}} {(n-1)\ld
r_{n-2}},~n\geq 2.\label{eq35}
\end{align}
\end{theorem}
\dproof By the orthogonality of $\phi\bl(\cdot,\cdot)$, we have
\begin{multline*}
\phi\bl(S\hl_{n+1}-\mu_nS\hl_n-\theta_{n-1}S\hl_{n-1},x^k) \\
=\ld k \la u_1, (P'_{n+1}-\tilde\sigma_nP'_n-\tilde\tau_{n-1}P'_{n-1})x^{k-1}\ra
=0, ~0\leq k\leq n-2,
\end{multline*}
or equivalently,
$$ \la u_1, (P'_{n+1}-\tilde\sigma_nP'_n-\tilde\tau_{n-1}P'_{n-1})R_k\ra=0,
~0\leq k\leq n-3,$$ so we have (\ref{eq33}) for some $\tilde d_n$ and
$\tilde e_n$, $n\geq 1$ with $\tilde d_0=\tilde e_0=\tilde e_{-1}=0$.
To find $\tilde e_n$, by using the orthogonality of
$\phi\bl(\cdot,\cdot)$, we have
\begin{align*}
&\phi\bl(S\hl_{n+1}-\mu_nS\hl_n-\theta_{n-1}S\hl_{n-1},x^{n-1}) \\
=&\la u_0,(P_{n+1}-\tilde\sigma_nP_n-\tilde\tau_{n-1}P_{n-1}) x^{n-1}\ra
\\
&+\ld(n-1)\la u_1,\{(n+1)R_n-\tilde d_{n}R_{n-1}-\tilde e_{n-1}R_{n-2}\} x^{n-2}\ra \\
=& -\tilde\tau_{n-1}p_{n-1}-\ld(n-1)\tilde e_{n-1}r_{n-2}
=-\theta_{n-1}s_{n-1}, ~n\geq 2
\end{align*}
so we have
$$\theta_{n-1}s_{n-1}=\tilde\tau_{n-1}p_{n-1}+(n-1)\ld \tilde e_{n-1}r_{n-2}, ~n\geq2.$$
This gives (\ref{eq35}).
To find $\tilde d_n$, by using the orthogonality of
$\phi\bl(\cdot,\cdot)$, we have, for $n\geq 1$,
\begin{equation}\label{eq36}
\begin{aligned}
&\phi\bl(S\hl_{n+1}-\mu_nS\hl_n-\theta_{n-1}S\hl_{n-1},S_n\hl)\\
=&\la u_0,(P_{n+1}-\tilde\sigma_nP_n-\tilde\tau_{n-1}P_{n-1})S\hl_n\ra \\
&+ \ld \la u_1, \{(n+1)R_n-\tilde d_{n}R_{n-1}-\tilde e_{n-1}R_{n-2}\} S^{(\ld)\prime}_n\ra= -\mu_n s_n .
\end{aligned}
\end{equation}
Since, from (\ref{eq31}),
$$ S\hl_n=P_n-\tilde\sigma_{n-1}P_{n-1}-\tilde\tau_{n-2}P_{n-2}+
\mu_{n-1}S\hl_{n-1}+\theta_{n-2}S\hl_{n-2},~n\geq 1,$$
we have the first part of (\ref{eq36}) as
\begin{multline}\label{eq37}
\la u_0, (P_{n+1}-\tilde\sigma_nP_n-\tilde\tau_{n-1}P_{n-1})S\hl_n)\ra \\
=-\tilde\sigma_np_n+\tilde\tau_{n-1}p_{n-1}(\tilde\sigma_{n-1}-\mu_{n-1}),~ n\geq 1.
\end{multline}
Since from (\ref{eq31}) and (\ref{eq33}),
$$ S^{(\ld)\prime}_n=(nR_{n-1}-\tilde d_{n-1}R_{n-2}-\tilde e_{n-2}R_{n-3})+
\mu_{n-1}S^{(\ld)\prime}_{n-1}+\theta_{n-2}S^{(\ld)\prime}_{n-2},~n\geq 1,$$
the second part of (\ref{eq36}) becomes
\begin{multline}\label{eq38}
\la u_1, \{(n+1)R_n-\tilde d_{n}R_{n-1}-\tilde e_{n-1}R_{n-2}\} S^{(\ld)\prime}_n\ra \\
= -n\tilde d_{n}r_{n-1}+\tilde d_{n-1}\tilde e_{n-1}r_{n-2}-(n-1)\tilde e_{n-1}\mu_{n-1} r_{n-2}, ~n\geq 1.
\end{multline}
Combining (\ref{eq37}) and (\ref{eq38}) with (\ref{eq36}), we have
\begin{multline*}
 \mu_ns_n =\tilde\sigma_np_n-(\tilde\sigma_{n-1}-\mu_{n-1})\tilde\tau_{n-1}p_{n-1} \\
+\ld \{n\tilde d_{n}r_{n-1}+((n-1)\mu_{n-1}-\tilde d_{n-1})\tilde e_{n-1}r_{n-2}\},~
n\geq 1.
\end{multline*}
Hence we have $\tilde d_{n}$, $n\geq 1$, in (\ref{eq34}). \qed

\vskip .5cm

Now we define extended generalized coherent pair.

\begin{definition}\label{def3.1}
Let $u_0$ and $u_1$ be quasi-definite
moment functionals with corresponding MOPS's $\pn$ and $\rn$,
respectively. Let $Q_n(x)=\frac{1}{n+1}P'_{n+1}(x),$ $n\geq 0.$

$\uu$ is an extended generalized coherent pair if there
exist complex numbers $d_n$, $e_n$, $\sigma_n,$ and $\tau_n,$ $n\geq0$ such that
\begin{multline}\label{eq3.9}
R_n(x) -d_{n-1}R_{n-1}(x)-e_{n-2}R_{n-2}(x) \\
=Q_{n}(x)-\sigma_{n-1} Q_{n-1}(x)-\tau_{n-2}Q_{n-2}(x),~n\geq 0, 
\end{multline}
where $d_{-1}=e_{-1}=e_{-2}=\sigma_{-1}=\tau_{-1}=\tau_{-2}=0$.

In particular, $\uu$ is an extended 3-term coherent pair if $e_n \tau_n\neq 0$, $n\geq 1$.
\end{definition}

Notice that if $e_n=0,$ $n\geq 1$, then extended generalized coherency is reduced to a simpler case such as extended coherency or $\rn=\qn$. 

This definition of an extended generalized coherent pair is equivalent to (\ref{eq33}).

\begin{remark}\label{rem32}
The following show that Theorem \ref{thm31} covers the
extended coherent pair and the extended symmetric coherent pair.
\begin{itemize}
    \item [(i)] If $\tilde\tau_n=0$ and $\theta_n=0$ for $n\geq 0$ in (\ref{eq31}), then $\uu$ is an extended coherent pair.
    \item[(ii)] Assume $u_0$ and $u_1$ are symmetric moment functionals. If  $\tilde\sigma_n=0$ and $\mu_n=0$ for $n\geq 0$ in (\ref{eq31}), then $\uu$ is an extended  symmetric coherent pair.
\end{itemize}
\end{remark}

\begin{theorem}\label{thm34}
$\uu$ is a generalized coherent pair if and only if the associated MOPS's $\sn$ and $\pn$
satisfy the relation in (\ref{eq31}) along with (\ref{eq32}).
\end{theorem}

\dproof Assume that the MOPS's $\sn$ and $\pn$ associated with the functionals $\uu$ satisfy (\ref{eq31}) and (\ref{eq32}). Then from (\ref{eq34}) and (\ref{eq35}), (\ref{eq32}) implies that $\tilde d_n=0$ and $\tilde e_n=0$ for $n\geq 1$ in (\ref{eq33}). The proof of the converse can be found in \cite{KLM}. \qed

From this theorem, notice that extended generalized coherency is an extension of generalized coherency.

\vskip .5cm

\section{Extended 3-term Coherent Pairs: The Classical Case}
\setcounter{equation}{0} \setcounter{theorem}{0}

In this section, let $\uu$ be an extended generalized coherent pair.
And let $\pn$ and $\rn$ be the MOPS's relative to quasi-definite
moment functionals $u_0$ and $u_1$, respectively, and satisfy extended 3-term coherency in (\ref{eq3.9}) with $\{Q_n(x):=\frac{P'_{n+1}(x)}{n+1}\}_{n=0}^\infty$.

We also assume that one of the moment functionals is classical, that is, the first case is when $u_0$ is classical, and the second case is when $u_1$ is strongly classical.

For the first case, consider that $u_0$ is classical with $(\varphi u_0)'=\psi u_0$. 
Then, it is well-known(\cite{Ch2}) that
$\{Q_n(x):=\frac{P'_{n+1}(x)}{n+1}\}_{n=0}^\infty$ is also a classical MOPS relative to $u=\varphi u_0$.

By acting $u$ and $u_1$ on (\ref{eq3.9}), we have
\begin{align}
&\left\{\begin{aligned} &\la u,R_n\ra= d_{n-1}\la
u,R_{n-1}\ra+ e_{n-2}\la u,
R_{n-2}\ra,~n\geq 3, \\
&\la u,R_2\ra= d_1\la u,R_1\ra+( e_0-\tau_0)\la
u,1\ra, \\
&\la u,R_1\ra=( d_0-\si_0)\la u,1\ra,
\end{aligned}\right. \label{eq52} \\
&\left\{\begin{aligned} &\la u_1,Q_n\ra=\si_{n-1}\la
u_1,Q_{n-1}\ra+\tau_{n-2}
\la u_1, Q_{n-2}\ra,~n\geq 3, \\
&\la u_1,Q_2\ra=\si_1\la u_1,Q_1\ra+(\tau_0 -
e_0)\la u_1,1\ra, \\
&\la u_1,Q_1\ra=(\si_0- d_0)\la u_1,1\ra.
\end{aligned}\right.\label{eq53}
\end{align}

\begin{theorem}\label{thm51}
Let $\uu$ be an extended generalized coherent pair, and let $u_0$
be a classical moment functional(or $u_1$ be a strongly classical one). Then
\begin{enumerate}
\item[(i)] if $\si_0= d_0$ and $\tau_0=
e_0$ (or equivalently, $\si_1=d_1$ and $\tau_1=e_1$), then $\rn
=\qn$ and $\si_n= d_n$ and $\tau_n= e_n$
for $n\geq 1$;

\item[(ii)] if $ e_n=\tau_n=0$, $n\geq 0$ (i.e.,
extended coherent pair) and $\si_0= d_0$, then $\rn
=\qn$ and $\si_n= d_n$;

\item[(iii)] if $\si_n= d_n=0$, $n\geq0$ (i.e.,
extended symmetric coherent pair) and $\tau_0= e_0$,
then $\rn =\qn$ and $\tau_n= e_n$ for $n\geq 1$.
\end{enumerate}
\end{theorem}

\dproof (i) From (\ref{eq52}) and (\ref{eq53}), the given
conditions $\si_0= d_0$ and $\tau_0= e_0$
give $\la u,R_n\ra=0$ and $\la u_1,Q_n\ra=0$, $n\geq 1$,
inductively. Thus $\rn=\qn$ and consequently $\si_n=
d_n$ and $\tau_n= e_n$ for $n\geq 1$. In the same way,
(ii) and (iii) (see \cite{KKMY}) hold. 
\qed

\vskip .5cm

In the following, we assume that $\uu$ is an extended 3-term coherent pair in (\ref{eq3.9}) with $e_n \tau_n \neq 0$, $n\geq 1$.
And we classify the 3-term extended coherent pairs when one of the moment functionals is classical or strongly classical. We find the coefficient relations between two corresponding MOPS's relative to the functionals.

We assume that $\qn$ and $\rn$ satisfy the following three-term
recurrence relation, for some $b_n, c_n, \beta_n, \gamma_n,$ $n\geq 0$,
\begin{align}
&xQ_n(x)=Q_{n+1}(x)+b_nQ_n(x)-c_nQ_{n-1}(x)~(c_n\neq0,~c_0=\la
u,1\ra),
\label{eq54}\\
&xR_n(x)=R_{n+1}(x)+\beta_nR_n(x)-\gamma_nR_{n-1}(x)~(\gamma_n \neq
0,~\gamma_0=\la u_1,1\ra). \label{eq55}
\end{align}

Now, we find the relations between coefficients, that is, $\{b_n, c_n\}$ for $\qn$, $\{\beta_n, \gamma_n \}$ for $\rn$, $\{\sigma_n, \tau_n\}$ and $\{d_n, e_n\}$ for the extended generalized coherent relation.

By multiplying (\ref{eq3.9}) by $x$, then using (\ref{eq54}) and
(\ref{eq55}), we have
\begin{equation}\label{eq56}
\begin{aligned}
&R_{n+1}+(\beta_n- d_{n-1})R_n+(\gamma_n-\beta_{n-1}
d_{n-1}- e_{n-2})R_{n-1}\\
&-( d_{n-1}+ e_{n-2}\beta_{n-2})R_{n-2}
- e_{n-2}\gamma_{n-2}R_{n-3} \\
=&Q_{n+1}+(b_n-\si_{n-1})Q_n+(c_n-b_{n-1}\si_{n-1}
-\tau_{n-2})Q_{n-1}\\
&-(\si_{n-1}+\tau_{n-2}b_{n-2})
Q_{n-2}-\tau_{n-2}c_{n-2}Q_{n-3},~n\geq 2.
\end{aligned}
\end{equation}

By applying (\ref{eq3.9}) to (\ref{eq56})
and then applying $R_{n-3}$ in (\ref{eq3.9}) with $n$ replaced by $n-1$, 
we have
\begin{equation}\label{eq58}
\begin{aligned}
&A_nR_n+B_nR_{n-1}-C_nR_{n-2} \\
&=D_nQ_n+E_nQ_{n-1}-F_nQ_{n-2}-G_nQ_{n-3},~n\geq 0,
\end{aligned}
\end{equation}
where
\begin{equation} \label{eq4.8}
\begin{aligned}
&\left\{ \begin{aligned} 
&A_n=\beta_n+ d_n-d_{n-1},\\
&D_n=b_n+\si_n-\si_{n-1},
\end{aligned}\right\}~n\geq 1, 
~ A_0=\beta_0+ d_0, ~D_0=b_0+\si_0, \\
&\left\{ \begin{aligned}
&B_n=\gamma_n-\beta_{n-1} d_{n-1}+  e_{n-1} -
e_{n-2}-\frac{ e_{n-2}\gamma_{n-2}}{ e_{n-3}},\\
&E_n=c_n-b_{n-1}\si_{n-1}+\tau_{n-1}-\tau_{n-2}
-\frac{ e_{n-2}\gamma_{n-2}}{ e_{n-3}},\\
\end{aligned}\right\}~n\geq 3,  \\
&~~~B_2=\gamma_2-\beta_1 d_1+ e_1- e_0,
~~~E_2=c_2-b_1\si_1 +\tau_1 -\tau_0,\\
&~~~B_1=\gamma_1-\beta_0 d_0 + e_0, ~E_1=c_1-b_0\si_0+\tau_0,\\
&\left\{ \begin{aligned}
&C_n= d_{n-1}+ e_{n-2} \beta_{n-2}-\frac{
e_{n-2}\gamma_{n-2}}{ e_{n-3}} d_{n-2},\\
&F_n=\si_{n-1}+\tau_{n-2}b_{n-2}-\frac{ e_{n-2}
\gamma_{n-2}}{ e_{n-3}}\si_{n-2},\\
 \end{aligned}\right\}~n\geq 3,  \\
&~~~C_2= d_1+ e_0 \beta_0,~F_2=\si_1+\tau_0 b_0, \\
&~~G_n=\tau_{n-2}c_{n-2}-\frac{ e_{n-2}\gamma_{n-2}}
{ e_{n-3}}\tau_{n-3},~n\geq 3,
\end{aligned}
\end{equation}
with $B_0=E_0=C_0=C_1=F_0=F_1=G_0=G_1=G_2=0$.

From (\ref{eq58}), since the leading coefficient of both sides are
the same, we have $A_n=D_n$, $n\geq 0$.

\vskip .3cm

{\bf Case I:} $A_n=D_n=0$, $n\geq 0$. \\
Then (\ref{eq58}) becomes
\begin{equation}\label{eq5.11}
B_nR_{n-1}-C_nR_{n-2}=E_nQ_{n-1}-F_nQ_{n-2}-G_nQ_{n-3},~ n\geq 1.
\end{equation}
Since the leading coefficients of both sides are the same, we have
$B_n=E_n$, $n\geq 1$.

{\bf Case I-1:} $B_n=E_n=0$, $n\geq 1$. \\
Then (\ref{eq5.11}) becomes
\begin{equation}\label{eq5.12}
C_nR_{n-2}=F_nQ_{n-2}+G_nQ_{n-3},~n\geq 2.
\end{equation}
Since the leading coefficients of both sides are the same, we have
$C_n=F_n$, $n\geq 2$.

{\bf Case I-1-(i):} $C_n=F_n=0$, $n\geq 2$. \\
Then $G_n=0$, $n\geq 3$. Hence we have from (\ref{eq4.8})
\begin{equation}\label{eq5.13}
\begin{aligned}
&A_n=D_n=0,~n\geq 0, ~~
&&B_n=E_n=0, ~n\geq 1,  \\ 
&C_n=F_n=0, ~n\geq 2,  
&&G_n=0,~n\geq 3. 
\end{aligned}
\end{equation}

\vskip 0.1cm 
{\bf Case I-1-(ii):} $C_n=F_n\neq 0$, $n\geq2$. \\
Then (\ref{eq5.12}) is to be
$$ R_{n-2}=Q_{n-2}+\frac{G_n}{C_n}Q_{n-3},~n\geq 2.$$
Hence with (\ref{eq3.9}), $G_n=0$, $n\geq 3$. This implies $\rn=\qn$.  

\vskip 0.1cm
{\bf Case I-2:} $B_n=E_n\neq 0$, $n\geq 1$. \\
Then (\ref{eq5.11}) is to be
$$ R_{n-1}-\frac{C_n}{B_n}R_{n-2}=Q_{n-1}-\frac{F_n}{B_n}Q_{n-2}
-\frac{G_n}{B_n}Q_{n-3},~n\geq 2.$$ Hence with (\ref{eq3.9}) this is reduced to trivial case or 2-term extended coherency, not the extended 3-term case.

\vskip .3cm

{\bf Case II:} $A_n=D_n\neq 0$, $n\geq 0$.\\
Then by (\ref{eq3.9}), (\ref{eq58}) is to be
\begin{multline}\label{eq5.24}
\left({B_n}+ d_{n-1}A_n\right)R_{n-1}-\left({C_n}-
e_{n-2}{A_n}\right)R_{n-2} \\
=\left({E_n}+\si_{n-1}A_n\right)Q_{n-1}
-\left({F_n}-\tau_{n-2}A_n\right)Q_{n-2}
-{G_n}Q_{n-3},~n\geq 1.
\end{multline}
Since the leading coefficients of both side are the same, we have
$$K_n:={B_n}+ d_{n-1}A_n
={E_n}+\si_{n-1}A_n,~n\geq 1.$$

{\bf Case II-1:} $K_n=0$, $n\geq 1$. \\
Then (\ref{eq5.24}) is to be
\begin{equation}\label{eq5.25}
\left({C_n}- e_{n-2}A_n\right)R_{n-2}
=\left({F_n}-\tau_{n-2}v\right)Q_{n-2}
+{G_n}Q_{n-3},~n\geq 2.
\end{equation}
Hence,
$$X_n:={C_n}- e_{n-2}A_n={F_n}-\tau_{n-2}A_n,~n\geq 2.$$

{\bf Case II-1-(i):} $X_n=0$, $n\geq 2$. \\ 
Then $G_n=0$, $n\geq 3$. Hence we have from (\ref{eq4.8})
\begin{equation}\label{eq5.26}
\begin{aligned}
&A_n=D_n(\neq0),~n\geq 0, \\
&K_n=B_n+d_{n-1}A_n=E_n+\si_{n-1}D_n=0, ~n\geq1\\ 
&X_n= C_n- e_{n-2}A_n=F_n-\tau_{n-1}D_n=0,~n\geq2,\\ 
& G_n=0,~n\geq 3.
\end{aligned}
\end{equation}

\vskip 0.1cm 
{\bf Case II-1-(ii):} $X_n\neq0$, $n\geq 2$. \\
Then (\ref{eq5.25}) is to be
$$ R_{n-2}=Q_{n-2}-\frac{G_n}{X_n} Q_{n-3},
~n\geq 2.$$ Then with (\ref{eq3.9}), $G_n=0$, $n\geq 3$. This implies $\rn=\qn$.

{\bf Case II-2:} $K_n\neq 0$, $n\geq 1$. \\
Then (\ref{eq5.24}) is to be
\begin{equation*} 
R_{n-1}-\frac{C_n- e_{n-2}}{K_n}R_{n-2} 
=Q_{n-1}-\frac{F_n-\tau_{n-2}}{K_n}Q_{n-2}
-\frac{G_n}{K_n}Q_{n-3},~n\geq 1.
\end{equation*}
Then with (\ref{eq3.9}), this is not the extended 3-term case.
\qed

\vskip .3cm

For the second case, we define $\tn$ as $R_n(x)=\frac{T'_{n+1}(x)}{n+1}$, then we integrate (\ref{eq3.9}). Hence we have
\begin{multline}
T_n(x)- d_{n-1}T_{n-1}(x)- e_{n-2}T_{n-2}(x)
\\
=P_n(x)-\si_{n-1} P_{n-1}(x)-\tau_{n-2}P_{n-2}(x),~n\geq 1,
\end{multline}
where $e_{-1}=\tau_{-1}=0$.
For $\tn$ to be classical, we assume $u_1$ is strongly classical. Then this case is the same as the first case with different notations. So we have the same result as the first case: when $u_0$ is classical.

So far, we find the relations concerned with the coefficients of extended generalized coherent pairs and the three-term recurrence relations of the corresponding MOPS when one of the moment functional is classical (or strongly classical). 

In this classification, there are four pairs of coefficients, that is, $\{d_n, e_n\}$ and $\{\sigma_n, \tau_n\}$ ($n\geq 0$) for extended 3-term coherent pairs in (\ref{eq3.9}) and $\{b_n, c_n\}$, $\{\beta_n, \gamma_n\}$ ($n\geq0$) for the three-term recurrence coefficients of $\qn$ and $\rn$, respectively in (\ref{eq54}) and (\ref{eq55}).
To find other coefficient relations for four pairs when some of the pairs are given, it is enough to give only one pair. 

In the above classification, there are two cases which are the extended 3-term coherency.

For the first case {\bf Case I-1-(i)}, (\ref{eq5.13}) gives
\begin{eqnarray}
 & \beta_n+ d_n-d_{n-1}=0,~n\geq1,~\beta_0+ d_0=0, \label{eq1-1} \\ 
 & b_n+\si_n-\si_{n-1}=0,~n\geq1,~ b_0+\si_0=0, \label{eq1-1b}
\end{eqnarray}
\begin{equation}\label{eq1-2}
\gamma_n-\beta_{n-1} d_{n-1}+  e_{n-1} -
e_{n-2}-\frac{ e_{n-2}\gamma_{n-2}}{ e_{n-3}}=0,~n\geq3, 
\end{equation}
$$\gamma_2-\beta_1 d_1+ e_1- e_0=\gamma_1-\beta_0 d_0 + e_0=0, $$
\begin{equation}\label{eq1-2c}
c_n-b_{n-1}\si_{n-1}+\tau_{n-1}-\tau_{n-2}
-\frac{ \tau_{n-2}c_{n-2}}{ \tau_{n-3}}=0,~n\geq3, 
\end{equation}
$$ c_2-b_1\si_1 +\tau_1 -\tau_0=c_1-b_0\si_0+\tau_0=0, $$
\begin{eqnarray}
&d_{n-1}+ e_{n-2} \beta_{n-2}-\frac{e_{n-2}\gamma_{n-2}}{ e_{n-3}} d_{n-2}=0, ~n\geq3, ~d_1+ e_0 \beta_0=0, \label{eq1-3} \\ 
&\si_{n-1}+\tau_{n-2}b_{n-2}-\frac{ e_{n-2} \gamma_{n-2}}{ e_{n-3}}\si_{n-2}=0, ~n\geq 3,~\si_1+\tau_0 b_0=0, \label{eq1-3s}
\end{eqnarray}
\begin{equation}\label{eq1-4}
\tau_{n-2}c_{n-2}e_{n-3}=e_{n-2}\gamma_{n-2}\tau_{n-3},~n\geq 3. 
\end{equation}

Now, we consider four cases that one of the four pairs is given, then we find coefficient relations for the other three pairs as follows.

{\bf Case 1 $\{\sigma_n, \tau_n\}$:} For the first case, when $\{\si_n, \tau_n\}$ are given, then from
(\ref{eq1-1b}) and (\ref{eq1-2c}), we have
\begin{align}
&~~~ b_n=\si_{n-1}-\si_n,~n\geq 0,\label{eq5.19}\\
&\left\{\begin{aligned}
&c_n=b_{n-1}\si_{n-1}\tau_{n-1}+\tau_{n-2}
+\frac{\tau_{n-2}c_{n-2}}{\tau_{n-3}},~n\geq 3,
\label{eq5.20} \\
& c_2=b_1\si_1-\tau_1+\tau_0,~c_1=b_0\si_0
-\tau_0.\end{aligned}\right.
\end{align}
Hence, we have $b_n$, $n\geq 0$ from (\ref{eq5.19}), and also we
have $c_n$, $n\geq 1$ from (\ref{eq5.19}) and (\ref{eq5.20}) recursively. 
Next, from (\ref{eq1-2}), (\ref{eq1-3}), and (\ref{eq1-4}), we have
\begin{align}
&\left\{ \begin{aligned}
& \gamma_n= e_{n-2}-
e_{n-1}+\beta_{n-1} d_{n-1}+\frac{
e_{n-2}\gamma_{n-2}}{ e_{n-3}},
~n\geq 3 \label{eq5.21} \\
&\gamma_2= e_0- e_1+\beta_0 d_1,
~\gamma_1=- e_0+\beta_0 d_0,\end{aligned}\right. \\
&\left\{\begin{aligned} &  d_{n-1}
=\frac{\tau_{n-2}c_{n-2}} {\tau_{n-3}}
d_{n-2}- e_{n-2} \beta_{n-2},
~n\geq 3, \label{eq5.22}\\
&  d_1=- e_0 \beta_0,\end{aligned}\right. \\
& ~~~ e_{n-2}=\frac{ e_{n-3}\tau_{n-2}}
{\gamma_{n-2}\tau_{n-3}}c_{n-2},~n\geq 3. \label{eq5.23}
\end{align}
Hence, from (\ref{eq1-1}) and (\ref{eq5.21})---(\ref{eq5.23}), we
find $\beta_n$, $\gamma_n$, $d_n$, and $e_n$, $n\geq1$, recursively with initial condition $ d_0$
and $ e_0$ (or equivalently $\beta_0$ and $\gamma_1$) in the
order of $[\beta_n, \gamma_{n+1},  e_{n+1}, d_{n+1}]$, $n\geq 1$.

{\bf Case 2 $\{b_n,c_n\}$:} Now, consider the second case by assuming that only $\{b_n, c_n\}$ are given, then (\ref{eq1-1b}) implies 
\begin{equation}\label{eqk1}
\sigma_n=\sigma_{n-1}-b_n,~n\geq 1,
\end{equation}
we get $\sigma_n$, $n\geq 1$ with initial condition $\si_0$. 
(\ref{eq1-2c}) and (\ref{eq1-4}) imply
\begin{equation}\label{eqk2}
\begin{aligned}
&\tau_{n-1}=\tau_{n-2}+b_n\si_{n-1}-c_n+
\frac{\tau_{n-2}c_{n-2}}{\tau_{n-3}},~n\geq 3, \\
& \tau_1=\tau_0+b_1\si_1-c_2,~
\tau_0=b_0\si_0-c_1, 
\end{aligned} 
\end{equation}
we have $\tau_n$, $n\geq 0$.
Then, in the same way as the first case, we get $\beta_n$,
$\gamma_n$, $d_n$, and $e_n$, $n\geq1$,
recursively with initial condition $ d_0$ and $ e_0$
(or equivalently $\beta_0$ and $\gamma_1$) in the order of
$[\beta_n,\gamma_{n+1}, e_{n+1},  d_{n+1}]$, $n\geq
1$.

{\bf Case 3 $\{\beta_n, \gamma_n\}$:} For third case, we assume that only $\{\beta_n$, $\gamma_n\}$ are given, then (\ref{eq1-1}) implies 
\begin{equation}\label{eqt1}
d_n=d_{n-1}-\beta_n,~n\geq 1,
\end{equation}
hence we get $d_n$, $n\geq 1$ with initial condition $d_0$. 
And (\ref{eq1-2}) implies
\begin{equation}\label{eqt2}
\begin{aligned}
&e_{n-1} =e_{n-2}-\gamma_n+\beta_{n-1} d_{n-1}+\frac{ e_{n-2}
\gamma_{n-2}}{ e_{n-3}}, ~n\geq3, \\
&e_1=e_0-\gamma_2 +\beta_1 d_1,
\end{aligned}
\end{equation}
hence we have $e_n$, $n\geq 1$ with the initial $e_0$.

Then, in the same way as the first case, we get $b_n$,
$c_n$, $\sigma_n$, and $\tau_n$, $n\geq1$,
recursively with initial condition $ d_0$ and $ e_0$
(or equivalently $\beta_0$ and $\gamma_1$) in the order of
$[b_n,c_{n+1}, \tau_{n+1},  \sigma_{n+1}]$, $n\geq1$.

{\bf Case 4 $\{d_n,e_n\}$:} Lastly, we assume $d_n, d_n$ are given, then we find from (\ref{eq1-1}) $\beta_n$, $n\geq 0$, from (\ref{eq1-2}) we find $\gamma_n$, $n\geq1$.
\begin{equation}\label{eq432}
\begin{aligned}
    & \beta_n=d_{n-1}-d_n, ~n\geq 1, ~\beta_0=- d_0 \\
    & \gamma_n=\beta_{n-1} d_{n-1}- e_{n-1} + e_{n-2} +\frac{ e_{n-2} \gamma_{n-2}}{ e_{n-3}},~ n\geq3,\\
    & \gamma_2=\beta_1d_1-e_1+e_0, ~\gamma_1=\beta_0 d_0-e_0.
    \end{aligned}
\end{equation}
Then, in the same way as the first case, we get $b_n$,
$c_n$, $\sigma_n$, and $\tau_n$, $n\geq1$,
recursively with initial condition $\sigma_0$ and $\tau_0$
in the order of $[b_n,c_{n+1}, \tau_{n+1},  \sigma_{n+1}]$, $n\geq1$.

For the second case {\bf Case II-1-(i)}, we have from (\ref{eq5.26})
\begin{equation}\label{eq2-1}
\beta_n+ d_n- d_{n-1} = b_n+\si_n-\si_{n-1}(\neq0),~n\geq 1,~\beta_0+ d_0 = b_0+\sigma_0,
\end{equation}
\begin{multline}\label{eq2-2}
    \gamma_n + e_{n-1}- e_{n-2} -\frac{ e_{n-2} \gamma_{n-2}}{e_{n-3}} \\
+ d_{n-1}(\beta_n-\beta_{n-1}+ d_n- d_{n-1} )=0, ~n\geq3,
\end{multline}
\begin{align*}
   & \gamma_2+ e_1-e_0 + d_1(\beta_2-\beta_1+ d_2- d_{1})=0, \\
   & \gamma_1 + e_0  +d_0(\beta_1-\beta_0+ d_1- d_{0})=0,
\end{align*}
\begin{multline} \label{eq2-2c}
c_n+\tau_{n-1}-\tau_{n-2}-\frac{ \tau_{n-2}c_{n-2}}{ \tau_{n-3}}  \\
 +\sigma_{n-1}(b_n-b_{n-1}+\sigma_n-\sigma_{n-1})=0,~ n\geq 3,
\end{multline}
\begin{align*}
&  c_2 +\tau_1 -\tau_0+\sigma_1(b_2-b_1+\sigma_2-\sigma_{1})=0, \\
& c_1+\tau_0+\sigma_0(b_1-b_0+\sigma_1-\sigma_{0})=0,
\end{align*}
\begin{equation}\label{eq2-3}
    d_{n-1}-\frac{e_{n-2}\gamma_{n-2}} { e_{n-3}} d_{n-2}- 
e_{n-2}(\beta_n- \beta_{n-2}+ d_n- d_{n-1})=0, ~n\geq 3,
\end{equation}
$$ d_1- e_0(\beta_2- \beta_0+ d_2- d_{1}) =0,$$
\begin{equation} \label{eq2-3s}
    \si_{n-1}-\frac{\tau_{n-2}
c_{n-2}}{ \tau_{n-3}}\sigma_{n-2}-\tau_{n-2}(b_n+\sigma_n-b_{n-2}-\sigma_{n-1})=0,~n\geq3,
\end{equation}
$$\si_1-\tau_0 (b_2- b_0+\sigma_2-\sigma_{1})=0,$$
\begin{equation}\label{eq2-4}
\tau_{n-2}c_{n-2}e_{n-3}= e_{n-2} \gamma_{n-2}\tau_{n-3}, ~n\geq 3.
\end{equation}
In this case, (\ref{eq2-1})---(\ref{eq2-4}) are related to each of the two pairs, (\ref{eq2-2}) and (\ref{eq2-3}) are related to the first two pairs $\{\beta_n, \gamma_n\}$ and $\{d_n, e_n\}$ and (\ref{eq2-2c}) and (\ref{eq2-3s}) are related to the second two pairs $\{b_n, c_n\}$ and $\{\sigma_n, \tau_n\}$.
Thus to solve four pairs of coefficients, we need at least two pairs of coefficients, of which one must be given from the first two pairs, and the other is given from the second two pairs.

From these classifications, we conclude with the following theorem. 
\begin{theorem}\label{thm50}
Let $v$ and $u$ be moment functionals with the corresponding MOPS's $\rn$ and $\qn$, respectively, are related by (\ref{eq3.9}) with $e_n \tau_n \neq 0$, $n\geq 0$. And $\qn$ and $\rn$ satisfy the three-term recurrence relations (\ref{eq54}) and (\ref{eq55}), respectively.

Then we have the relation (\ref{eq58}), $A_n=D_n$, $n\geq 0$, and $G_n=0$, $n\geq3$, that is, 
$$  \beta_n+ d_n- d_{n-1}  =b_n+\si_n-\si_{n-1}, ~n\geq0.$$
$$e_{n-3}\tau_{n-2}c_{n-2} = e_{n-2}  \tau_{n-3}\gamma_{n-2},~n\geq 3.$$

Moreover, there are only two cases for extended 3-term coherent pairs.
\begin{itemize}
    \item[(I)] When $A_n=D_n=0,$ $n\geq0$, $B_n=E_n=0,$ $n\geq1$, and $C_n=F_n=0,$ $n\geq2$ in (\ref{eq58}), the coefficients for the above four pairs have the relations (\ref{eq1-1})---(\ref{eq1-4}).    
    
    And for getting four pairs of coefficients, it is enough that only one of the pairs is given, then we find other three pairs of coefficients relations as {\bf Case 1 --- Case 4}
    \item[(II)] When $A_n=D_n\neq0,$ $n\geq0$, ${B_n}+A_n d_{n-1} ={E_n}+A_n\si_{n-1}=0,$ $n\geq1$, and ${C_n}-A_n e_{n-2} ={F_n} -A_n \tau_{n-2} =0, ~n\geq 2,$
    the coefficients for the above four pairs have the relations (\ref{eq2-1})---(\ref{eq2-4}).

    And for getting four pairs of coefficients, it is enough that only two of the pairs are given, of which one is one of $\{\beta_n, \gamma_n\}$ and $\{d_n, e_n\}$ and the other is one of $\{b_n, c_n\}$ and $\{\sigma_n, \tau_n\}$, then we find other two pairs of coefficients relations from (\ref{eq2-1})---(\ref{eq2-4}).
\end{itemize}
\end{theorem}

\vskip .5cm

Now, we study the companion moment functional of the extended 3-term coherent pairs when $u_0$ is classical or when $u_1$ is strongly classical. When $u_1$ is strong classical, we have the same result as when $u_0$ is classical. Hence we state the following:

\begin{theorem}\label{thm4.3}
Let $v$ and $u$ be moment functionals with the corresponding MOPS's $\rn$ and $\qn$, respectively, related by (\ref{eq3.9}). Then there exist polynomials $A(x)$ and $D(x)$ of degree 2, satisfying
\begin{equation}\label{eqk3}
A(x) v=D(x) u.
\end{equation}
\end{theorem}

\dproof
Let $A(x)$ and $D(x)$ be polynomials of degree 2.
We define $M_k:=\la A(x)v, R_k(x)\ra$, $k\geq 0$.

By using (\ref{eq3.9}), we have
\begin{equation}\label{eq6.1}
\la A(x)v,Q_n(x)\ra=\sigma_{n-1}\la A(x)v,Q_{n-1}\ra
+\tau_{n-2}\la A(x)v, Q_{n-2}\ra, ~n\geq 5
\end{equation}
and
\begin{equation}\label{eq6.2}
\begin{aligned}
{\rm (i)}~ \la A(x)v, Q_0(x)\ra =&\la A(x)v,R_0(x)\ra=M_0,\\
{\rm (ii)}~ \la A(x)v, Q_1(x)\ra =&M_1+(\sigma_0-{d}_0)M_0,\\
{\rm (iii)}~\la A(x)v, Q_2(x)\ra =&M_2+(\sigma_1-{d}_1)M_1+ \{ \sigma_1 (\sigma_0 -{d}_0)+(\tau_0-{e}_0)\}M_0,\\
{\rm (iv)}~\la A(x)v, Q_3(x)\ra =& M_3+ (\sigma_2-{d}_2)M_2+\{ \sigma_2 (\sigma_1-{d}_1)+(\tau_1-{e}_1)\}M_1 \\
&+\{\sigma_2(\sigma_1(\sigma_0-{d}_0)+( \tau_0-{e}_0))+\tau_1(\sigma_0-{d}_0)\}M_0,\\
{\rm(v)}~\la A(x)v, Q_4(x)\ra =& M_4-{d}_3M_3+(\tau_2-{e}_2)M_2 + \tau_2(\sigma_1-{d}_1)M_1 \\
&+\tau_2 \{ \sigma_1 (\sigma_0-{d}_0)+(\tau_0-{e}_0)\}M_0.
\end{aligned}
\end{equation}
Since $M_k=0$, $k\geq 3$, only $M_0$, $M_1$, and $M_2$ remain.

To satisfy (\ref{eq6.1}) inductively, $\la A(x)v, Q_3(x)\ra =\la A(x)v, Q_4(x)\ra =0$ hold in (\ref{eq6.2}).

Since $\deg(A(x))=2$, we have $M_2=\la A(x) v, R_2(x)\ra \neq 0$. Hence we can solve $M_0$ and $M_1$ from (\ref{eq6.2}) (iv) and (v). Then these $M_0$ and $M_1$ make (\ref{eq6.2}) (iii) $\frac{{e}_2M_2}{\tau_2}\neq 0$.
Thus, we have a polynomial $D(x)$ of degree 2, satisfying $A(x)v=D(x)u$ and (\ref{eq6.2}).
\qed

\vskip .3cm

\begin{remark}\label{rem4.2}
(i) When a strongly classical moment functional $v$ with MOPS $\qn$ and a classical moment functional $u_0$ with MOPS $\pn$ satisfy $Q_n(x)=\frac{1}{n+1}P'_{n+1}(x)$, then $(\varphi u)'=\psi u$ and $\varphi u_0 =u$. Hence Theorem \ref{thm4.3} implies, if $\uu$ is an extended 3-term coherent pair, then 
\begin{equation}\label{eq4.33}
A(x) u_1=D(x) u=D(x) \varphi u_0,
\end{equation}
and $2\leq\deg(D(x) \varphi) \leq 4$, in detail
$$ \deg(D(x) \varphi)=\left\{ \begin{aligned}
&2, ~{\rm if~} u_0=u_H, \\
&3, ~{\rm if~} u_0=u_L, \\
&4, ~{\rm if~} u_0=u_J, u_B. 
\end{aligned} \right.
$$
(ii) When $u_0$ and $u_1$ are a symmetric coherent pair, then Theorem \ref{thm4.3} gives the result of Theorem 3.5 in \cite{DM2}.
\end{remark}

\vskip .5cm

When only $\{b_n, c_n\}$ are given in {\bf Case 2}, we describe the algorithm in order to obtain $\sigma_n$, $\tau_n$, $\beta_n$, $\gamma_n$, $d_n$, and $e_n$, $n\geq1$. Note that $\sigma_n$ and $\tau_n$, $n\geq1$, can be obtained without concerning other parameters. So we construct the algorithm with two parts.

\begin{enumerate}
\item[{STEP I}] FROM (\ref{eqk1}) FINDING $\{\sigma_n\}$ AND FROM (\ref{eqk2}) FINDING$\{\tau_n\}$

\begin{itemize}
\item[{ }] STARTING INITIAL VALUES:
$$ \sigma_0$$
\item[{STEP I-1}]
FROM (\ref{eqk1}) WE GET $\sigma_n$  RECURSIVELY BY
$$ \sigma_n=\sigma_{n-1}-b_n$$
\item[{STEP I-2}]
FROM (\ref{eqk2}) WE GET $\tau_n$ RECURSIVELY BY
$$\tau_0=b_0\si_0-c_1$$
$$ \tau_1=\tau_0+b_1\si_1-c_2$$
$$\tau_{n-1}=\tau_{n-2}+b_n\si_{n-1}-c_n+\frac{\tau_{n-2}c_{n-2}}{\tau_{n-3}}$$

\end{itemize}

\item[{STEP II}] FINDING $\{\beta_n\}$, $\{\gamma_n\}$, $\{ d_n\}$, AND $\{ e_n\}$

\begin{itemize}
\item[{ }] STARTING INITIAL VALUES:
$$ \beta_0 \qquad \gamma_1$$
\item[{STEP II-1}]
FROM (\ref{eq1-1}) AND (\ref{eq5.21}) FINDING
$$ d_0=-\beta_0 \qquad  e_0=-\gamma_1 -\beta_0^2$$
\item[{STEP II-2}]
FROM (\ref{eq5.23}) AND (\ref{eq5.22}), WE FIND THE FOLLOWING ONE BY ONE
$$  e_{1} \qquad  d_{1} $$
\item[{STEP II-3}]
FROM (\ref{eq5.13}), (\ref{eq5.21}), (\ref{eq5.23}), AND (\ref{eq5.22}), WE FIND THE FOLLOWING ONE BY ONE
$$ \beta_1 \qquad \gamma_{2} \qquad  e_{2} \qquad  d_{2} $$
$$ \begin{matrix}  \vdots \end{matrix} $$
\item[{STEP II-$n$}]
FROM (\ref{eq5.13}), (\ref{eq5.21}), (\ref{eq5.23}), AND (\ref{eq5.22}), WE FIND THE FOLLOWING ONE BY ONE
$$ \beta_n \qquad \gamma_{n+1} \qquad  e_{n+1} \qquad  d_{n+1} $$
\end{itemize}
\end{enumerate}
\vskip .5cm

\begin{example}\label{eg4-1}
Let $u_0=u_J^{(1,0)}$ be a Jacobi moment functional, that is,
$$ \la u_0, p(x)\ra :=\int_{-1}^1 p(x) (1-x)dx.~\forall p(x)\in \mathbb{P}.$$
Then $u=\varphi u_0=(1-x^2)u_0$ and the corresponding OPS is $\{P^{(2,1)}_n(x)\}_{n=0}^\infty$. So from \cite{Ra}, the MOPS $\qn$ relative to $u$ is
$$Q_n(x)=\frac{n!(n+3)!}{(2n+3)!} \sum_{k=0}^n\binom{n+2}{n-k} \binom{n+1}{k} (x-1)^k(x+1)^{n-k},~n\geq 0,$$
and $\qn$ satisfy the following three-term recurrence relation
$$ x Q_n(x) =Q_{n+1}(x) +b_n Q_n(x) -c_n Q_{n-1}(x), ~n\geq 0,$$
$$ b_n= \frac{-3}{(2n+3)(2n+5)},~~ c_n=\frac{-(n+1)(n+3)}{(2n+3)^2}.$$

\begin{itemize}
    \item [(i)] To find coefficients for four pairs, this is {\bf Case 2: $\{b_n,c_n\}$}. From STEP I of the above algorithm and (\ref{eqk1}) and (\ref{eqk2}), 
    we can find the parameters $\sigma_n$ and $\tau_n$ with initial $\sigma_0$ of recurrence relation in (\ref{eq3.9}) as the following
    $$ \sigma_n= \sigma_0+\frac{3}{2(2n+3)} -\frac{3}{14},~n\geq 1,$$
    \begin{multline*} 
    \tau_n = \frac{-3}{(2 n + 5) (2 n + 7)}\si_0+ \frac{4n^2 (28 n + 207) + 804n + 669}{28 (2 n + 3) (2 n + 5)^2 (2 n + 7)} \\
    + \tau_{n-1} -\frac{ n (n + 2)}{(2 n + 1)^2 }\frac{\tau_{n - 1}}{\tau_{n - 2}} + \frac{1}{4},~ n\geq 2, 
    \end{multline*}
    $$\tau_1=\frac{-2\sigma_0}{7}-\frac{306}{1225}, ~~\tau_0=\frac{-\sigma_0}{5}+\frac{8}{25}.$$
    Also, we can find from STEP II of the above algorithm $\beta_n$ and $\gamma_n$ which are the three-term recurrence relation coefficients of $\rn$ in (\ref{eq54}) and $d_n$ and $e_n$ in (\ref{eq3.9}) the coefficients of extended 3-term coherent pairs, recursively.
    \item[(ii)] We assume $A(x)=(1-x)^2$ and $D(x)=(1+x)^2$ in (\ref{eqk3}) with $v=u_1$ under the same conditions as in (i), then
    $$ u_1=\frac{(1+x)^2}{(1-x)^2}u+M\delta(1-x)+N\delta'(1-x)
    =u_J^{(0,3)}+M\delta(1-x)+N\delta'(1-x),$$
    where $M=\langle u_1,1\rangle=\gamma_0$ and $N=\langle u_1, 1-x\rangle=\gamma_0(1-\beta_0)$ are constants.(see Lemma2.4(iii) in \cite{KLM} and \cite{LK})
    Therefore, $u_1$ is a perturbed Jacobi moment functional.
    
    \item[(iii)] We assume $A(x)=(1-x)(1+x)$ and $D(x)=(1-x)(1+x)$ in (\ref{eqk3}) with $v=u_1$ under the same conditions as in (i), then
    $$ u_1=u+M\delta(1-x)+N\delta(1+x)
    =u_J^{(2,1)}+M\delta(1-x)+N\delta(1+x),$$
    where $M=\frac{1}{2}\langle u_1,1-x\rangle=\frac{1}{2}\gamma_0(1-\beta_0)$ and $N=-\frac{1}{2}\langle u_1, 1+x\rangle=-\frac{1}{2}\gamma_0(1+\beta_0)$ are constants.(see Lemma2.4(ii) in \cite{KLM} and \cite{LK})
    Therefore, $u_1$ is a perturbed Jacobi moment functional(see \cite{Ko}).
\end{itemize}
We have different companion moment functionals in (ii) and (iii). This is because of different initial values $\sigma_0$.
\qed
\end{example}

\begin{example}\label{eq4-2}
Let $u_1=u_L^{(2)}$ be a Laguerre moment functional, that is,
$$ \la u_0, p(x)\ra :=\int_{0}^\infty p(x) x^2e^{-x}dx.~\forall p(x)\in \mathbb{P}.$$
Then the corresponding OPS is $\{L^{(2)}_n(x)\}_{n=0}^\infty$. So from \cite{Ra}, the MOPS $\rn$ relative to $u_1$ is
$$R_n(x)=  \sum_{k=0}^n(-1)^{n-k} \frac{n!}{k!} \binom{n+2}{n-k} x^k ,~n\geq 0,$$
and $\qn$ satisfy the following three-term recurrence relation
$$ x R_n(x) =R_{n+1}(x) +\beta_n R_n(x) -\gamma_n R_{n-1}(x), ~n\geq 0,$$
$$ \beta_n= 2n+3,~~ \gamma_n=n(n+2).$$
\begin{itemize}
    \item [(i)] To find the coefficients for four pairs, this case is {\bf Case 3: $\{\beta_n, \gamma_n\}$}. From (\ref{eq5.21})---(\ref{eq5.23}),
    we can find the parameters $d_n$ and $e_n$ with initial $d_0$ and $e_0$ as the following
    $$ d_n= d_0+n(n+4),~n\geq 1,$$
    $$ e_{n-1}=e_{n-2}-\frac{e_{n-2}}{e_{n-3}} n(n-2)+d_0(2n+3)+2n(n^2+5n+5), ~n\geq3,$$
    $$ e_1 = e_0+5d_0+17.$$
    
    Also, we can find from STEP II of the above algorithm $\beta_n$ and $\gamma_n$ which are the three-term recurrence relation coefficients of $\rn$ in (\ref{eq54}) and $\sigma_n$ and $\tau_n$ in (\ref{eq3.9}) the coefficients of extended 3-term coherent pairs, recursively.
    \item[(ii)] We assume $A(x)=x^2$ and $D(x)=(1+x)^2$ in (\ref{eqk3}) with $v=u_1$ under the same conditions as in (i), then
    \begin{multline*}
    u=\frac{x^2}{(1+x)^2}u_1+M\delta(1+x)+N\delta'(1+x)
    \\=\frac{u_L^{(4)}}{(1+x)^2}+M\delta(1+x)+N\delta'(1+x),
    \end{multline*}
    where $M=\langle u,1\rangle=c_0$ and $N=\langle u, 1+x\rangle=c_0(1+b_0)$ are constants.(see Lemma2.4(iii) in \cite{KLM} and \cite{LK})
    Moreover, the companion of $u_1$ is $\phi u_0=u$, that is $u_0=\frac{1}{x}u+K\delta(x)$, where $K=\langle u_0, 1\rangle$.(see Lemma2.4(i) in \cite{KLM} and \cite{LK})
    
    Therefore, $u$ and $u_0$ are different from the perturbed Laguerre moment functionals.
\qed
\end{itemize}
\end{example}

\end{document}